\documentclass[conference]{IEEEtran}
\ifCLASSINFOpdf
\else
\fi

\usepackage{mathrsfs}
\usepackage{graphicx}
\usepackage{eurosym}
\usepackage{amssymb}
\usepackage{amsmath}
\usepackage{amsfonts}
\usepackage{epstopdf}
\usepackage{epsf,subfigure}
\usepackage{psfrag}
\usepackage{graphics}
\usepackage{color} 
\usepackage{cite}
\usepackage{slashbox}
\usepackage{array,multirow,pbox}
\usepackage{enumitem}

\newcommand{\rem}[1]{}

\newtheorem{definition}{Definition}
\newtheorem{example}{Example}

\newcommand{\qed}{\nobreak \ifvmode \relax \else
      \ifdim\lastskip<1.5em \hskip-\lastskip
      \hskip1.5em plus0em minus0.5em \fi \nobreak
      \vrule height0.75em width0.5em depth0.25em\fi}

\graphicspath{{Figures/}}

\newcommand{\mysubeq}[2]{
\begin{subequations}\label{#1}
\begin{align}
#2
\end{align}\end{subequations}}

\begin{document}
%
\title{Prioritized Threshold Allocation for Distributed Frequency Response}



%
\author{\IEEEauthorblockN{Sai Pushpak Nandanoori\IEEEauthorrefmark{1}\IEEEauthorrefmark{2},
Soumya Kundu\IEEEauthorrefmark{2},
Draguna Vrabie\IEEEauthorrefmark{2},
Karan Kalsi\IEEEauthorrefmark{2} and
Jianming Lian\IEEEauthorrefmark{2}}
\IEEEauthorblockA{\IEEEauthorrefmark{1}Electrical \& Computer Engineering, Iowa State University, Ames, IA 50011 USA\\ Email: pushpak@iastate.edu}
\IEEEauthorblockA{\IEEEauthorrefmark{2}Optimization \& Control Group, Pacific Northwest National Laboratory, Richland, WA 99354 USA\\
Email:\,\{soumya.kundu,\,draguna.vrabie,\,karanjit.kalsi,\,jianming.lian\}@pnnl.gov}}


\maketitle

\begin{abstract}
Higher penetration of renewable generation will increase the demand for adequate (and cost-effective) controllable resources on the grid that can mitigate and contain the contingencies locally before it can cause a network-wide collapse. However, end-use constraints can potentially lead to load unavailability when an event occurs, leading to unreliable demand response services. Sensors measurements and knowledge of the local load dynamics could be leveraged to improve the performance of load control algorithms. In the context of hierarchical frequency response using ensemble of switching loads, we present a metric to evaluate the fitness of each device in successfully providing the ancillary service. Furthermore a fitness-based assignment of control set-points is formulated which achieves reliable performance under different operating conditions. Monte Carlo simulations of ensembles of electric water heaters and residential air-conditioners are performed to evaluate the proposed control algorithm.
\end{abstract}


%
\IEEEpeerreviewmaketitle

\section{Introduction}

In today's power systems operations, traditional frequency control resources (e.g. speed governors, spinning reserves) are deployed to ensure resilient grid operations under contingencies, by restoring system frequency close to its nominal values. The importance of adequate (and cost-effective) frequency response mechanisms is expected to grow even further as the grid turns `greener' and `smarter'. Electrical loads, if coordinated smartly, have the potential to provide a much faster, cleaner and less expensive alternative to the traditional frequency responsive resources \cite{Strbac:2002,Short:2007,Callaway:2011}. Potential of controllable loads to provide frequency response services has been explored both in the academia and the industry \cite{Schweppe:1982,Kirby:2002,Callaway:2011, Kundu:2011PSCC, Perfumo:2012, Sinitsyn:2013, Mathieu:2013, Ma:13, Zhang:2013, Hao:14, Sanandaji:16}.

In order to scalably intergate millions of controllable devices (loads) into the grid operational paradigm, a hierarchical distributed control architecture is conceptualized in which a supervisor (e.g. a load aggregator) is tasked with dispersing the response of the loads across the ensemble so that some desirable collective behavior is attained. Dispersion of load response in frequency (by assigning to each load specific frequency thresholds to respond to \cite{Lu:2006,Horst:2007,Molina_Garcia:2011,Lian:2016,Kundu:2017Assessment}) allows certain power-frequency droop-like response enabling easier integration of such frequency-responsive resources in the grid operational framework. However, the availability of the end-use load to respond to ancillary service requests strongly depends on the local dynamics and constraints, adversely affecting the reliability of decentralized frequency control algorithms. 

In this article, we consider a hierarchical control framework for coordinating ensembles of switching loads (electric water-heaters and residential air-conditioners) to provide frequency response services to the grid operator. An aggregator is tasked with coordinating the devices' response to provide frequency response when an event happens. Each device receives a frequency threshold which it uses to turn `on' or `off' autonomously when a frequency event happens, if it is permitted to do so by the end-use level constraints. However, due to the extremely short duration and unpredictability of the frequency event, it is impractical to coordinate and communicate with the devices in real-time \textit{after} the event occurrence. Instead, the aggregator has to estimate the availability of devices that can turn `off' or `on' on-demand without violating end-use constraints. This paper proposes a metric to evaluate the \textit{`fitness'} of each device in the aggregation to respond to frequency response requests over any given period. The `fitness' values are then used by the aggregator to assign thresholds in a prioritized way such that when a frequency event happens, the `fittest' devices are likeliest to respond first. Sec.\,\ref{S:problem} describes the systems and the control problem. Sec.\,\ref{S:fitness} introduces the `fitness' metric, while Sec.\,\ref{S:priority} discusses fitness-based prioritization of loads for frequency response. In Sec.\,\ref{S:results} we present simulation results, before concluding the article in Sec.\,\ref{S:concl}.

\section{System Description}\label{S:problem}

\subsection{Load Models}

Consider an ensemble of $N$ switching devices which continuously switch between two operational states to serve certain local demand. Air-conditioners and electric water-heaters fall under such categories. Such devices can be modeled in general as,
\mysubeq{E:x}{
\dot{x}^i(t)&=a^i\,x^i(t) + b^i(t) + c^i\,p^i(t)\,,~\forall t\geq0\,, \\
p^i(t^+)&=\left\lbrace \begin{array}{cl}
0\,, & \text{if }h_1^ix^i(t)+h_2^i(t)\geq\delta^i/2\\
P^i\,, & \text{if }h_1^ix^i(t)+h_2^i(t)\leq-\delta^i/2\\
p^i(t)\,, & \text{otherwise}
\end{array}\right..\label{E:s}
}
where $x^i(t)$ is the vector of the dynamic state variables; $p^i(t)\!\in\!\lbrace 0,P^i\rbrace$ represents the discrete power consumption of the device with power rating $P^i$\,; $a^i,\,b^i(t),\,c^i,\,h_1^i,\,h_2^i(t)$ and $\delta^i$ represent various device parameters and exogenous inputs (such as weather conditions). In particular, $\delta^i$ is often referred to as the hysteresis bandwidth. The power consumption of the device toggles between two values ($0$ in the `off' state and $P^i$ in the `on' state) according to some state-dependent switching condition as defined in \eqref{E:s}. Specific examples concerning residential air-conditioners and electric water heaters will be provided later in this article. The total power consumed by the ensemble is 
\begin{align}
p_\Sigma^{}(t):=\sum_{i=1}^N p^i(t)\,,~0\leq p_\Sigma^{}(t)\leq \sum_{i=1}^NP^i\,.
\end{align}
While $\sum_{i=1}^NP^i$ and $0$ denote, respectively, the maximal and minimal power that the population of devices can draw at any given moment, the actual flexibility refers to quantifying {the amount of change in behavior, from nominal, the ensemble of devices is capable of accommodating without compromising on end-user quality of service.}

\subsection{Hierarchical Control Architecture}

\begin{figure*}[thpb]
  \begin{center}
\includegraphics[scale=0.25]{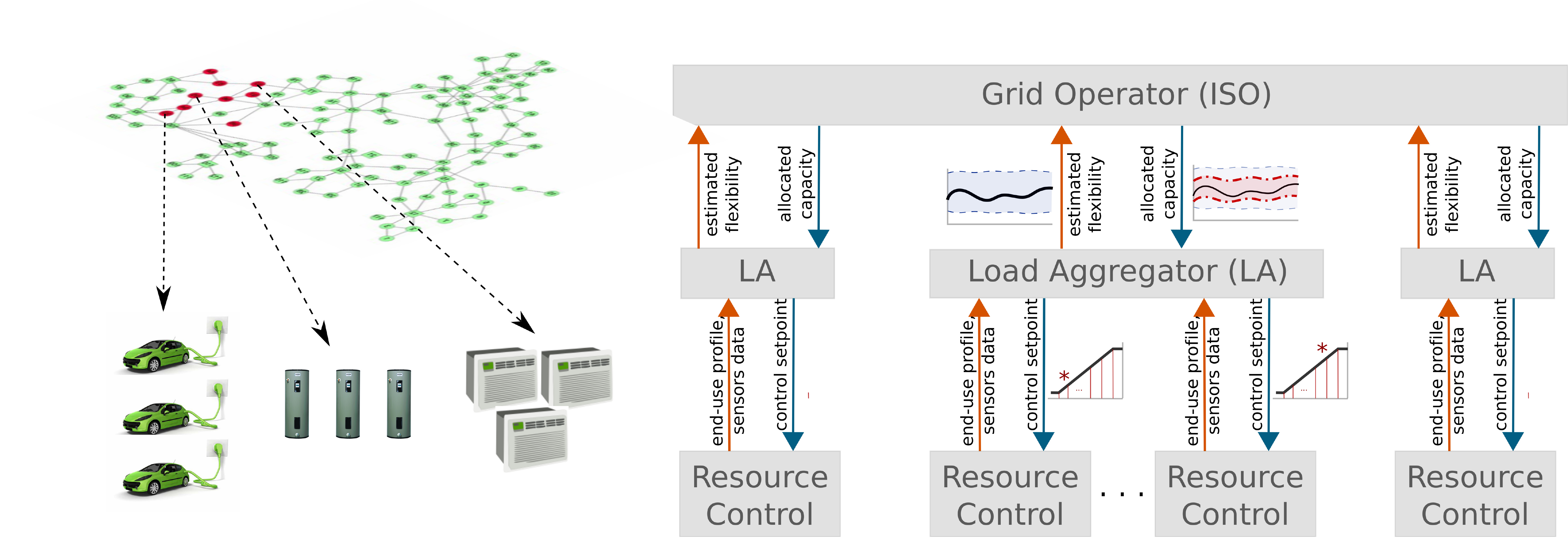} 
  \end{center}
  \caption{Hierarchical control architecture for engaging flexible end-use devices in grid frequency response services, via load aggregators.}
  \label{F:layers}
\end{figure*} 

The outline of the hierarchical control architecture is presented in Fig.\,\ref{F:layers}. The participation of controllable resources into the grid ancillary services at the Independent System Operator (ISO) layer is coordinated by \textit{load aggregators} \cite{Callaway:2011,Kundu:2017Assessment}. Load aggregators serve as the communication and control link between the \textit{resource controller} at the device-level and the \textit{grid operator}. It is the responsibility of the load aggregators to characterize the aggregated flexibility of the resources and implement control strategies to coordinate the pool of resources to provide the requested grid ancillary service, subject to customer quality-of-service (QoS) constraints. The flexibility could be represented in the form of the total amount by which the group of devices under an aggregator can increase and decrease their consumption over a given period of time, a portion of the reported flexibility is then procured/dispatched by the grid operator for frequency responsive services. At the start of each control window (5-15 min), new sensor measurements are used to update the control setpoints for each load across the grid, in a scalable manner. The aggregators received certain information (measurements from local sensors, predicted end-use profiles) from the devices and report to the grid operator the estimated flexibility in the aggregated power consumption over the control window. The grid operator analyzes the available flexibility across the network, along with the current network status (generation and load forecast, measurements from PMUs, topology) to optimally allocate the control capacities from the responsive load ensembles. Load aggregators then communicate back to the individual loads their updated control setpoints which the loads respond to autonomously, during the control window.

\begin{figure}
  \begin{center}
\includegraphics[scale=0.3]{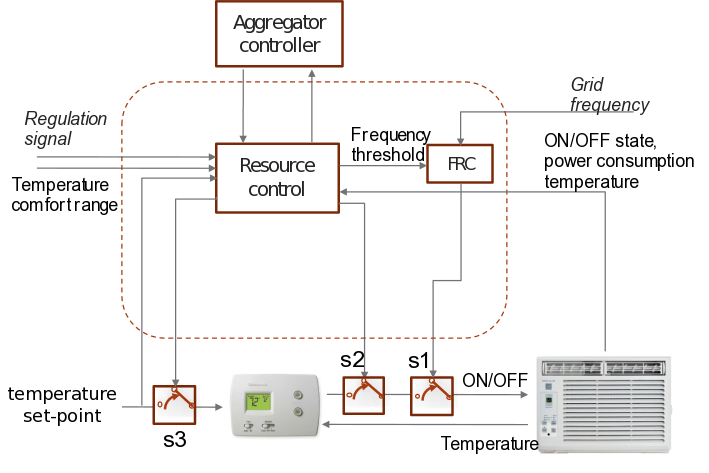} 
  \end{center}
  \caption{Example of a resource controller integrated with a thermostat-based control loop and equipped with a frequency response controller (FRC).}
  \label{fig:ACResourceControl}
\end{figure} 

The \textit{resource controller} acts as the interface between the existing {(local) device controller} and the load aggregator. It is the responsibility of the resource controller to characterize the flexibility and response capability of the device, and augment the existing device control to enable grid service response by executing the control setpoints received from the aggregator. Fig.\,\ref{fig:ACResourceControl} briefly depicts the control and communication pathways between an aggregator, resource controller and the device-level controller. In this diagram, the resource controller is augmented with the device-level thermostatic control, which switched the device on/off based on some temperature setpoint, while also having a \textit{frequency response controller} module that senses the grid frequency and based on some pre-assigned frequency threshold decides to either switch on/off the device, only if that forced switching does not violate the local temperature constraints (i.e. the end-use constraints take precedence over any grid service requests).

\subsection{Autonomous Frequency Response}

\begin{figure}[thpb]
\centering
\includegraphics[scale=0.3]{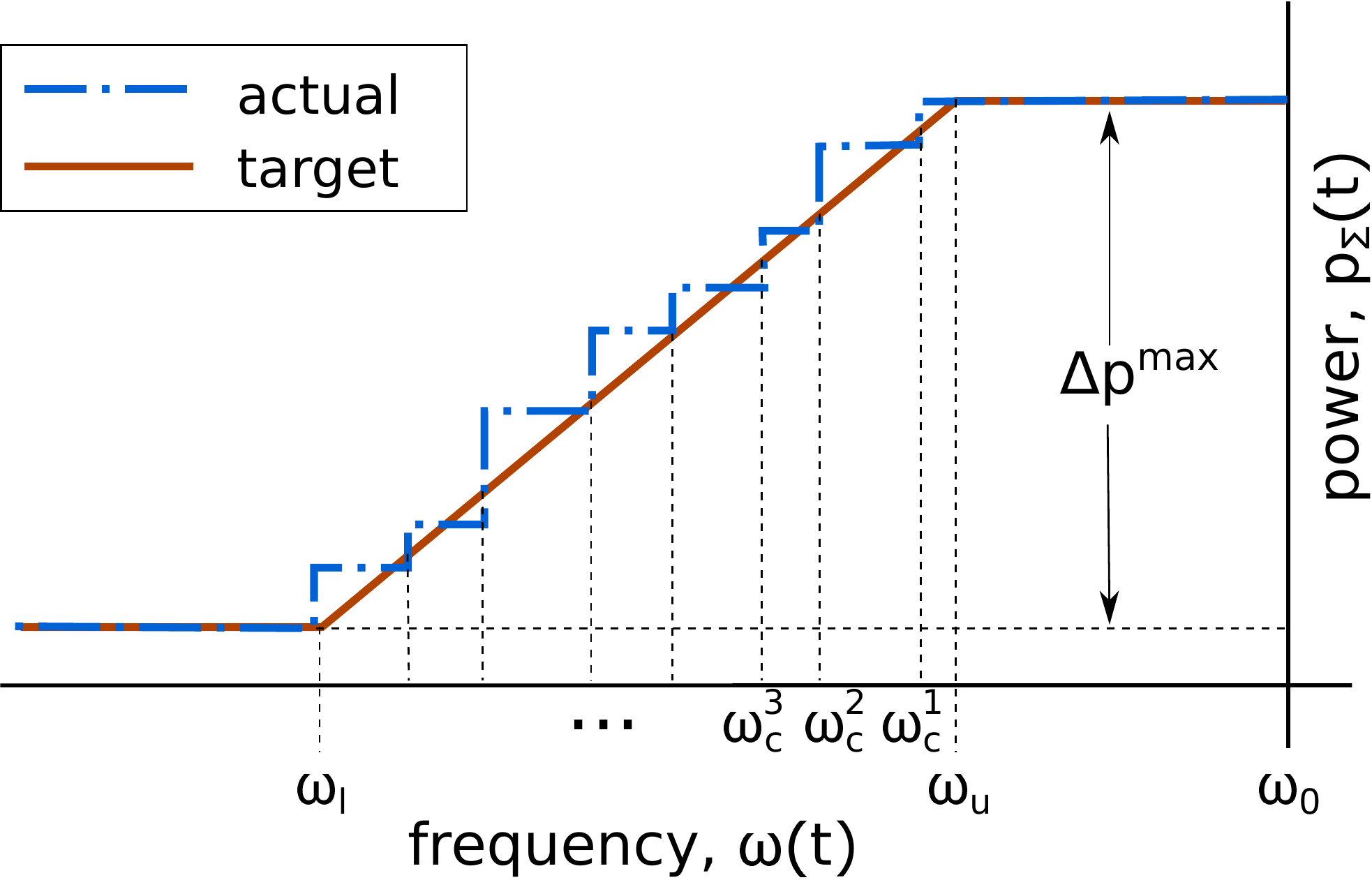}
\caption{Illustration of a power-frequency response curve.}
\label{F:droop}
\end{figure}

In this section, we explain the decentralized frequency response algorithm executed by the resource controller module. While similar approach can also be used for frequency regulation problem, the focus of this discussion will be only on the frequency response. Consider an ensemble of $N$ devices. When this ensemble commits to an under-frequency response, it is expected to decrease its power consumption by turning off some of its devices if the frequency falls. Since there are other frequency control mechanisms in-place, any such ensemble of loads will be expected to respond to events when the frequency is in a specific range. Thus a typical under-frequency response curve would look like Fig.\,\ref{F:droop}, where $\omega_u$ and $\omega_l$ denote the upper and lower limits of the frequency range assigned to the ensemble, and $\omega_0$ is the nominal frequency (60\,Hz). 
Clearly, $\omega_l<\omega_u\leq\omega_0$ (for under-frequency response). 
The target frequency response curve is a smooth line whose slope is determined based on the number (and power consumption) of the devices available to switch their states from `on' to `off'. The actual control is implemented by assigning frequency thresholds to each device, such that each device can turn `off' by monitoring the frequency on its own (see \cite{Molina_Garcia:2011,Lian:2016,Kundu:2017Assessment} for details). An over-frequency response policy can be constructed in a similar way. 

In order to better explain the response policy, let us assume that any given time $t$\,, $\mathcal{S}_t=\lbrace d_1,\,d_2,\dots,\,d_m\rbrace$ represents the set of indices of the devices that have committed to under-frequency response service. The target response capacity, i.e. the height of the power-frequency response curve in Fig.\,\ref{F:droop} is given by,
\begin{align}\label{E:capacity}
\Delta p^{\max}&:=\sum_{i=1}^mP^{d_i}\,.
\end{align}
Without any loss of generality, let us assume that the corresponding frequency thresholds $\left\lbrace \omega^i_c(t)\right\rbrace_{i=1}^m$ are chosen in an ordered way so that,
\begin{align}\label{E:w_priority}
\omega_l\leq\omega_c^m(t)<\dots<\omega_c^2(t)<\omega_c^1(t)\leq\omega_u.
\end{align} 
One possible way to choose the frequency thresholds to produce the target response curve in Fig.\,\ref{F:droop} is to assign 
\begin{align*}
\forall i\in\lbrace 1,\dots,m\rbrace:~\omega_c^i(t):=\omega_u-\frac{\omega_u-\omega_l}{\Delta p^{\max}}\sum_{j=1}^iP^{d_j}\,.
\end{align*}
The available devices obey the following response policy:
\begin{align}\label{E:response}
p^i(t^+)&=0\,,~\text{ if }\left\lbrace \begin{array}{cc}\omega(t)\leq\omega_c^i(t)\,\text{ and }\\
h_1^ix^i(t)+h_2^i(t)>-\delta^i/2\end{array}\right.
\end{align}
where the switching condition $h_1^ix^i(t)+h_2^i(t)>-\delta^i/2$ ensures that local end-use constraint does not require the device to be `on' at the time of response. If the device was indeed in a switched `on' state just before the response, $p^i(t)=P^i$\,, then by turning `off' the device successfully performs frequency response task. We have ignored, for simplicity, any finite time-delay in the response. The total power consumption of the ensemble under this response policy, as illustrated in Fig.\,\ref{F:droop}, is given by,
\begin{align}
\forall t:~p_\Sigma^{}(t^+)&=p_\Sigma^{}(t)-{\sum}_{\lbrace\forall d_i\in\mathcal{S}_t:\,\omega_c^i(t)\geq\omega(t)\rbrace}P^{d_i}\,.
\end{align}

\subsection{Problem Statement}

The key point to note here is that the values of the frequency thresholds in \eqref{E:w_priority} depend on the availability of the devices to turn `off' during an under-frequency event (equivalently, turn `on' for over-frequency event). However, continuous monitoring of the device states in an ensemble has high telemetry requirements, along with potential privacy concerns (for the device owners). A more viable option is to acquire and update the device states information once (at the start of) every fixed control time window, while using that information to estimate the availability of the responsive devices during the control window. 

This brings us to the problem of quantifying the \textit{fitness} of each of device to provide a certain kind of service (frequency response in this case). The fitness values can be used for the prioritized assignment of the frequency thresholds. For example, the devices that have high \textit{fitness} value (defined later) will be assigned the thresholds that are closer to the \textit{nominal} frequency, so that when the frequency starts to deviate the `fittest' devices are called for service first, while the thresholds that are farther away are assigned to the devices with lower fitness values. The focus of this article is to propose a metric to quantify the `fitness' of each controllable device to provide autonomous frequency response. Moreover, we will also explore how the fitness values can be used to estimate the maximum available flexibility in an ensemble of responsive devices to provide frequency response service.

\section{Metric for Fitness}\label{S:fitness}
\textit{Fitness} is a qualitative measure of the ability of a controllable load to (successfully) respond to a certain kind of ancillary service request. 

\begin{definition}\label{D:fitness}
`Fitness' of a device-$i$ for a particular service request-$k$ is denoted by a scalar $\pi^i_k\in[0,1]$ which quantifies how likely the device is to successfully complete the service request over some time window. 
\end{definition}

Evaluation of this metric of \textit{fitness} could depend on several factors, such as the device dynamics, response delays, rate of failure (to respond) and the type of service request. In the simplest of its form, the fitness could be composed of two component metrics - 1) \textit{availability} for response, and 2) \textit{quality} of response. 

\begin{definition}\label{D:availability}
`Availability' metric of a device-$i$ for a particular service request-$k$\,, denoted by $\pi^{\text{avail},i}_k\in[0,1]$\,, quantifies the probability that the device is available to respond to the particular service request over some time window. 
\end{definition}

\begin{definition}\label{D:quality}
`Quality' metric of a device-$i$ for a particular service request-$k$\,, denoted by $\pi^{\text{qual},i}_k\in[0,1]$\,, quantifies the probability that the device, when available, completes the service request successfully. 
\end{definition}

Consider for example, a plug-in electric vehicle that had signed up for a frequency response service at $t=t_0$ over a time window $[t_0\,,\,t_f]\,.$ Suddenly, at some time $t=t_0+\Delta t<t_f$ the vehicle was taken off the charger and driven away. If an event happens any time between $t_0+\Delta t$ and $t_f$\,, the vehicle would be unavailable to respond to it, even though it signed up for it. On the other hand, consider a residential air-conditioner that signs up for frequency response service over some control window. When the event happens, it is ready to respond to it. However, due to delays in control it takes a while to actually switch its mode of operation, thereby delivering a poor quality of service even though it was available to respond to it. The \textit{fitness} metric can be considered as a product of the \textit{availability} and \textit{quality} metrics as follows:
\begin{align}
\pi_k^i&=\pi^{\text{avail},i}_k\cdot\pi^{\text{qual},i}_k\,.
\end{align}
The success of a device in completing a service request, meeting all the performance metrics, depends on a lot of factors, such as sensors and actuation time-delays, as well as sensors and actuators failures. It is expected that during real-life implementations, the performance of a device to a particular service request could be monitored over time to estimate the \textit{quality} metric, $\pi^{\text{qual},i}_k$\,. As an example, in a simple form, the performance degradation due to time-delays can be modeled into the \textit{quality} metric as,
\begin{align}
\pi^{\text{qual},i}_k = \exp\left(-\beta\,t_{d,k}^i\right)
\end{align}
where $\beta>0$ is an appropriate scaling factor and $t_{d,k}^i\geq 0$ is an estimated (possibly over previous such requests) time-delay of the device in responding to the particular service request. A total failure would be captured by the limiting case $t_{d,k}^i\rightarrow+\infty$ while $t_{d,k}^i=0$ would refer to a success rate of 1. Both the availability and quality metrics can be estimated and updated online by monitoring the device performance in response to similar requests. The fitness values thus estimated can be maintained online at the aggregator level, and updated at the start of every control window and used in frequency threshold allocation. In the rest of this section, we illustrate how, in presence of sufficient device-level information, the \textit{availability} metrics can be computed for under- and over-frequency response.

\subsection{Under-Frequency Response}
When an under-frequency event happens, the devices that have previously committed for under-frequency response service are requested to turn `off' one-by-one based on their frequency thresholds. Thus the availability of a device to respond to an under-frequency response can be quantified by the probability that the device is in the `on' state when the event happens. It can be argued that the conditional probability distribution of the time of occurrence of a frequency event given the event has occurred within some interval $[t_0,t_f]$, is uniform over the time interval. Then the \textit{availability} metric of a device with respect to under-frequency response (denoted by the subscript `resp-') can be obtained as:
\begin{align}\label{E:avail_resp-}
\pi_{\text{resp-}}^{\text{avail},i}\!&=\!\mathbf{Pr}\lbrace \text{device-$i$ is `on' when under-frequency happens}\rbrace\nonumber\\
&=\!\int_{t_0}^{t_F}\mathbf{Pr}\lbrace \text{device-$i$ is `on' at time $\tau$}\rbrace\cdot f_{\text{resp-}}(\tau)\,d\tau\nonumber\\
&=\!\int_{t_0}^{t_f}\frac{s^i(\tau)}{t_f-t_0}\,d\tau = \frac{t^i_{on}}{t_f-t_0}
\end{align}
where, $s^i(\cdot)\in\lbrace 0,1\rbrace$ represents the operational state of the device, taking value $0$ in the `off' state and $1$ in the `on' state; $t_{on}^i$ is the length of time the device spends in the `on' state during the control window $[t_0,t_f]$\,; and $f_{\text{resp-}}(\cdot)$ is the uniform conditional probability density function of the time of occurrence of the under-frequency event, given that the event occurs during the interval. Of course, it is not possible to exactly know the `on' time length $t_{on}^i$ for each device. But with the knowledge of the internal states of the devices, and some forecast of the external conditions, it is possible to estimate $\pi_{\text{resp-}}^{\text{avail},i}$ for each device at the start of each control period.

\subsection{Over-Frequency Response}

Using similar arguments, the availability factor for an over-frequency response (denoted by the subscript `resp+') over an control window $[t_0,t_f]$ could be calculated as:
\begin{align}\label{E:avail_resp+}
\pi_{\text{resp+}}^{\text{avail},i}&= \frac{t^i_{off}}{t_f-t_0}
\end{align}
where $t_{off}^i$ is the length of time the device spends in the `off' state during the control window $[t_0,t_f]$\,.

\begin{example} (Air-conditioners) Consider an ensemble of $N$ air-conditioning (AC) loads. Each device dynamics is represented by \cite{Perfumo:2012},
\mysubeq{E:x_AC}{
\dot{T}(t)&=-\frac{\left(T(t)-T_a(t)\right)}{C\,R}-\frac{\eta\,p(t)}{C}\,,\\
p(t^+)&=\left\lbrace \begin{array}{cl}
0\,, & \text{if }T(t)\leq T_{set}-\delta T/2\\
P\,, & \text{if }T(t)\geq T_{set}+\delta T/2\\
p(t)\,, & \text{otherwise}
\end{array}\right.,
}
where $T(t)$ is the room temperature; $p(t)\in\lbrace 0,P\rbrace$ represent the power draw of the AC; $T_a(t)$ denotes the outside air temperature; and $C,\,R,\,\eta$ are the device parameters representing the room thermal resistance, thermal capacitance and the load efficiency, respectively. $T_{set}$ is the temperature set-point and $\delta T$ represents the width of the temperature hysteresis deadband. 
If the outside air temperature is constant throughout the control window, the dynamics \eqref{E:x_AC} can be solved to compute the time $t_{on}^{off}$ an initially `on' device spends before turning `off', and the time $t_{off}^{on}$ an initially `off' device spends before turning `on' as:
\begin{subequations}
\begin{align}
t_{on}^{off}&=CR\log\left[\frac{\left(T(t_0)-T_a(t)\right)+\eta PR}{\left(T_{set}-\delta T/2-T_a(t)\right)+\eta PR}\right],\\
t_{off}^{on}&=CR\log\left[\frac{T(t_0)-T_a(t)}{T_{set}+\delta T/2-T_a(t)}\right].
\end{align}
\end{subequations}
The time a device spends in the `on' and `off' state during the control window is given by:
\begin{subequations}
\begin{align}
t_{on}&=\left\lbrace \begin{array}{ll} \min\left(t_f\!-\!t_0\,,\,t_{on}^{off}\right)\,,&\text{`on' at $t_0$}\\
\max\left(0\,,\,t_f\!-\!t_0\!-\!t_{off}^{on}\right)\,,&\text{`off' at $t_0$}
\end{array}\right.\\
t_{off}&=t_f-t_0-t_{on}\,.
\end{align}
\end{subequations}
\hfill\hfill\qed\end{example}

\begin{example} (Electric water-heaters) Consider an ensemble of $N$ electric water-heating (EWH) loads. The water temerature dynamics of an EWH can be modeled using a `one-mass' thermal model which assumes that the temperature inside the water-tank is spatially uniform (valid when the tank is \textit{nearly} full or \textit{nearly} empty) \cite{Diao:2012}:
\begin{align}\label{E:x_EWH}
\dot{T}_w(t) &=-a(t)\,T_w(t)+b(s(t),t)\,,\\
\text{where, }a(t) &:=\frac{1}{C_w}\left(\dot{m}(t)\,C_p+W\right),\notag\\
\&~b(s(t),t)&:=\frac{1}{C_w}\left(s(t)\,P+\dot{m}(t)\,C_p\,T_{in}(t)+W\,T_a(t)\right).\notag
\end{align}
$T_w$ denotes the temperature of the water in the tank, and $s(t)$ denotes a switching variable which determines whether the EWH is drawing power ($s(t)=1$ or `on') or not ($s(t)=0$ or `off'). 
The state of the EWH (`on' or `off') is determined by the switching condition:
\begin{align}\label{E:s_EWH}
s(t^+)&=\left\lbrace \begin{array}{cl}
0\,, & \text{if }T_w(t)\geq T_{set}+\delta T/2\\
1\,, & \text{if }T_w(t)\leq T_{set}-\delta T/2\\
s(t)\,, & \text{otherwise}
\end{array}\right.,
\end{align}
where $T_{set}$ is the temperature set-point of the EWH with a deadband width of $\delta T$\,. 
If the exogenous parameters are unchanged throughout the control window, then the dynamics \eqref{E:x_EWH} can be solved to compute the time $t_{on}^{off}$ an initially `on' device spends before turning `off', and the time $t_{off}^{on}$ an initially `off' device spends before turning `on' as:
\begin{subequations}
\begin{align}
t_{on}^{off}&=\frac{1}{a(t)}\log\left[\frac{-a(t)\,T_w(t_0)+b(1,t)}{-a(t)\left(T_{set}+\delta T/2\right)+b(1,t)}\right],\\
t_{off}^{on}&=\frac{1}{a(t)}\log\left[\frac{-a(t)\,T_w(t_0)+b(0,t)}{-a(t)\left(T_{set}-\delta T/2\right)+b(0,t)}\right].
\end{align}\end{subequations}
the time the device spends in the `on' and `off' state during the control window as:
\begin{subequations}
\begin{align}
t_{on}&=\left\lbrace \begin{array}{ll} \min\left(t_f\!-\!t_0\,,\,t_{on}^{off}\right)\,,&\text{`on' at $t_0$}\\
\max\left(0\,,\,t_f\!-\!t_0\!-\!t_{off}^{on}\right)\,,&\text{`off' at $t_0$}
\end{array}\right.\\
t_{off}&=t_f-t_0-t_{on}\,.
\end{align}
\end{subequations}
 
\hfill\hfill\qed\end{example}

Note that, under perfect information, the availability metric can be computed exactly using the `on' time and `off' time duration. In more realistic scenarios the equations above can be used to estimate certain statistical properties (such as mean, variance) of the random variables $t_{on}$ and $t_{off}$ from any information on statistics of the sensor errors and device parameters.

\section{Prioritized Frequency Response}\label{S:priority}
At the start of the control window ($t=t_0$), the \textit{fitness} values of each device in the population are computed for any particular service-$k$\,. Based on the fitness values, $\pi_k^i\,\forall i\in\lbrace 1,2,\dots,N\rbrace$\,, all the devices in the population are prioritized in an order $\lbrace d_1,\,d_2,\,d_3,\,\dots,\,d_N\rbrace$ for consideration of commitment to service-$k$, such that 
\begin{align}\label{E:priority}
\pi^{d_1}_k\geq \pi^{d_2}_k\geq \pi^{d_3}_k\dots\geq \pi^{d_N}_k\,.
\end{align}
Accordingly the frequency thresholds (for primary frequency response) are assigned, such that thresholds closer to the nominal frequency are assigned to higher priority (`fitter') devices. Finally, a subset of those devices are selected, based on the priority order, such that their aggregate power rating equals (within some tolerable error, $\epsilon_p$) the target response capacity $\Delta p^{\max}$, i.e. choose the smallest $m\in\lbrace 1,2,\dots,N\rbrace$ such that:
\begin{align}\label{E:threshold}
\left|\Delta p^{\max}-\sum_{i=1}^m P^{d_i}\right|\leq\epsilon_p\,.
\end{align}
Furthermore, the probability that all the devices succeed in responding to the service request is given by,
\begin{align}
\mathbf{Pr}\lbrace \text{ `success' }\rbrace=\prod_{i=1}^m\pi_k^{d_i}\,,
\end{align}
while probability that at least one device fails to respond to the request is given by
\begin{align}
\mathbf{Pr}\lbrace \text{ at least one `failure' }\rbrace\geq\left(1-\pi_k^{d_m}\right)\,,
\end{align}
since $\pi_k^{d_1}\geq\pi_k^{d_2}\geq\dots\geq\pi_k^{d_m}$\,.
Consequently, we can also compute the maximal response capacity that the aggregation of devices can commit to with a probability of success 1, which is defined as follows.
\begin{definition}
The maximal capacity that an aggregation can successfully (i.e. with probability 1) commit for any service-$k$ is denoted by $\left[\Delta p\right]_k^{cap}$ and is given by
\begin{align}\label{E:maxcap}
\left[\Delta p\right]_k^{cap}&:=P^{d_1}+P^{d_2}+\dots+P^{d_n}\,,
\end{align}
such that $\pi_k^{d_i}=1\,\forall i\in\lbrace 1,2,\dots,n\rbrace$ and $\pi_k^{d_{n+1}}<1$\,.
\end{definition}
\begin{figure}[thpb]
\centering
\includegraphics[scale=0.2]{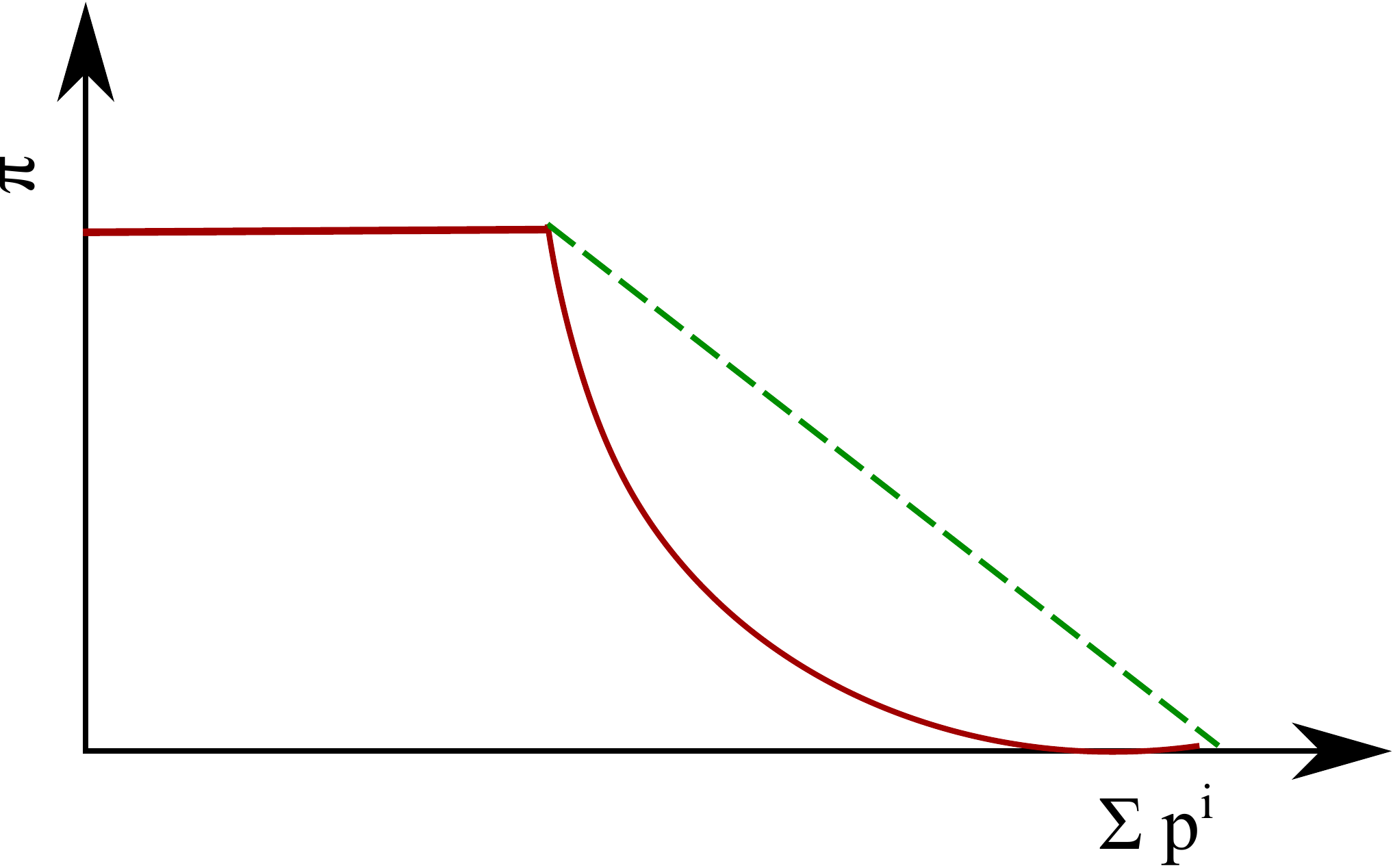}
\caption{Illustration of availability metric as a function of load power.}
\label{F:example}
\end{figure}
Moreover the availability metric as a function of the aggregated load power (as illustrated in Fig.\,\ref{F:example}) could signify how sensitive the performance is to the accuracy of the estimated maximal flexibility (up to $\pi_k^i=1$). The steeper slope (`red' curve) signifies that the performance drops faster (as opposed to the `green dashed' line) the estimated available flexibility. Future work will address this issue more closely.
\subsection{Control Performance}

Control performance will be evaluated against a metric termed as the \textit{reserve margin variability target} (RMVT) which is expressed as the following,
\begin{align*}
RMVT&:=\left|1-\frac{\text{total response provided on request}}{\text{total response requested}}\right|\,.
\end{align*}
There are several sources of uncertainties that may affect the control performance (and contribute to the RMVT), such as the forecast errors, modeling uncertainties and faults in sensors. Of particular interest to this work are the uncertainties due to control parameters, such as the discrete allocation of frequency thresholds and sampling time delays. 

\subsubsection{Discrete Threshold Allocation}
Note from Fig.\,\ref{F:droop} that, due to the discrete power consumption of the devices and due to finite number of devices committing to provide the frequency response, at any frequency value, the maximal deviation is given by the power rating of the largest device that has signed up for the service. Using the threshold assignment logic given in \eqref{E:priority} and \eqref{E:threshold}, the relative error in response (relative to the response capacity) due to discrete loads can be given by,
\begin{align}
\frac{\delta p}{\Delta p^{\max}}\leq \frac{\max_{i=1,\dots,m}P^{d_i}}{\sum_{i=1}^mP^{d_i}}
\end{align}
This suggests that in order for the relative error in response to be less than certain pre-specified performance metric, $\epsilon$, the committed aggregate capacity for response has to be larger than certain critical value, i.e.
\begin{align}\label{E:err_discrete}
{RMVT}\leq \epsilon\implies \Delta p^{\max}\geq \frac{\max_iP^i}{\epsilon}
\end{align}

\subsubsection{Sampling Time}
\begin{figure}[thpb]
\centering
\includegraphics[scale=0.3]{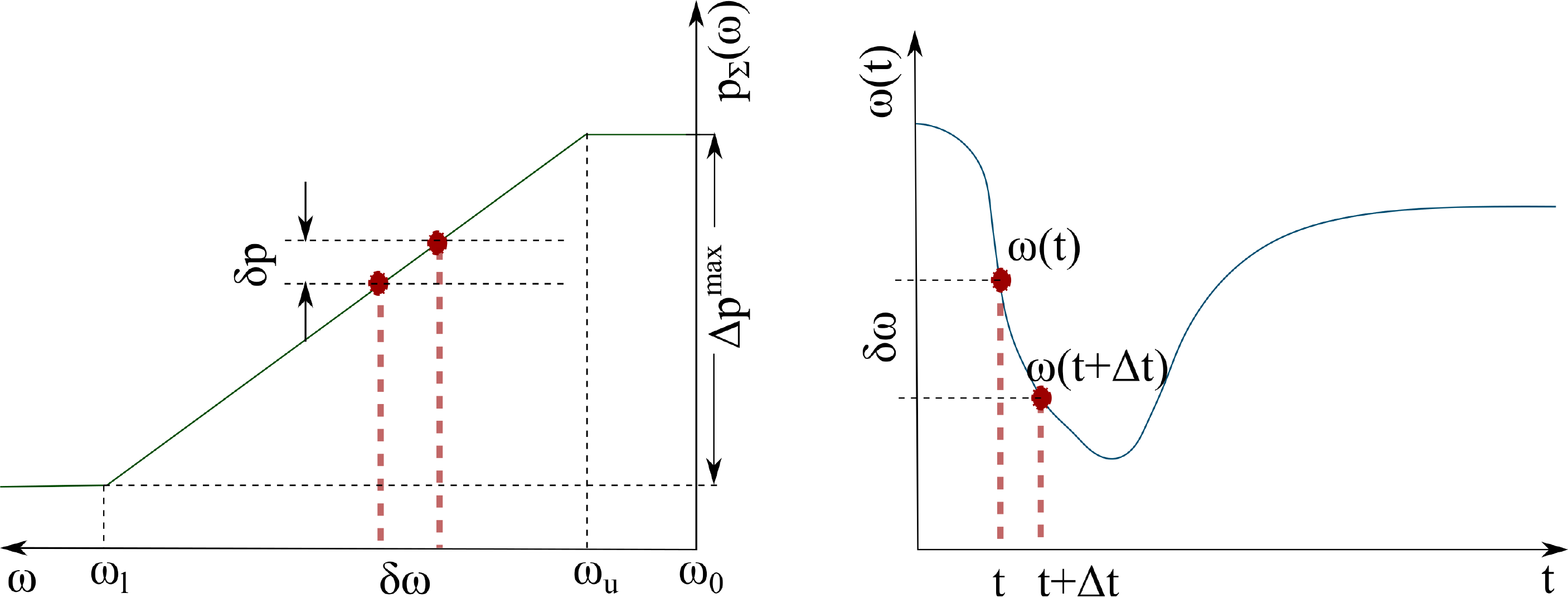}
\caption{Response error due to finite non-zero sampling time.}
\label{F:timedelay}
\end{figure}
In general, the delay between a sensor measurement or service request and the control execution (or, actuation) can have multiple components, e.g. actuator delays ($t_d$), and delays due to finite sampling rate ($\Delta t$). The effective delay in response, therefore will be given by $\max\left(t_d\,,\Delta t\right)$. For simplicity, in this work, the actuator delay is assumed to be negligible ($t_d=0$), while we focus only on the delay due to finite sampling rate. Consider the event shown in Fig.\,\ref{F:timedelay}. The device local sensor detects the that frequency $\omega(t)$ at time $t$ is lower than its assigned threshold and therefore, according to the control logic \eqref{E:response}, turns `off' at the next time instant $t+\Delta t$\,, where $\Delta t$ refers to the non-zero discrete sampling time. If the frequency deviates further during this time delay, there is a non-zero deviation, $\delta p$\,, from the target frequency response curve. For small $\Delta t$\,, this error can be estimated at,
\begin{subequations}\label{E:err_sampling}
\begin{align}
\delta p \approx \left.\frac{\partial p_\Sigma^{}(\omega)}{\partial \omega}\right|_{\omega(t)}\cdot\left.\frac{\partial \omega}{\partial t}\right|_t\cdot\Delta t\\
\implies \frac{\delta p}{\Delta p^{max}}\approx \frac{\Delta t}{(\omega_u-\omega_l)}\left.\frac{\partial \omega}{\partial t}\right|_t
\end{align}\end{subequations}
i.e. the relative error in response due to finite sampling rate is proportional to the sampling time, as well as the rate of change of frequency. 

\section{Simulation Results}\label{S:results}
\begin{figure}[thpb]
\centering
\includegraphics[scale=0.48]{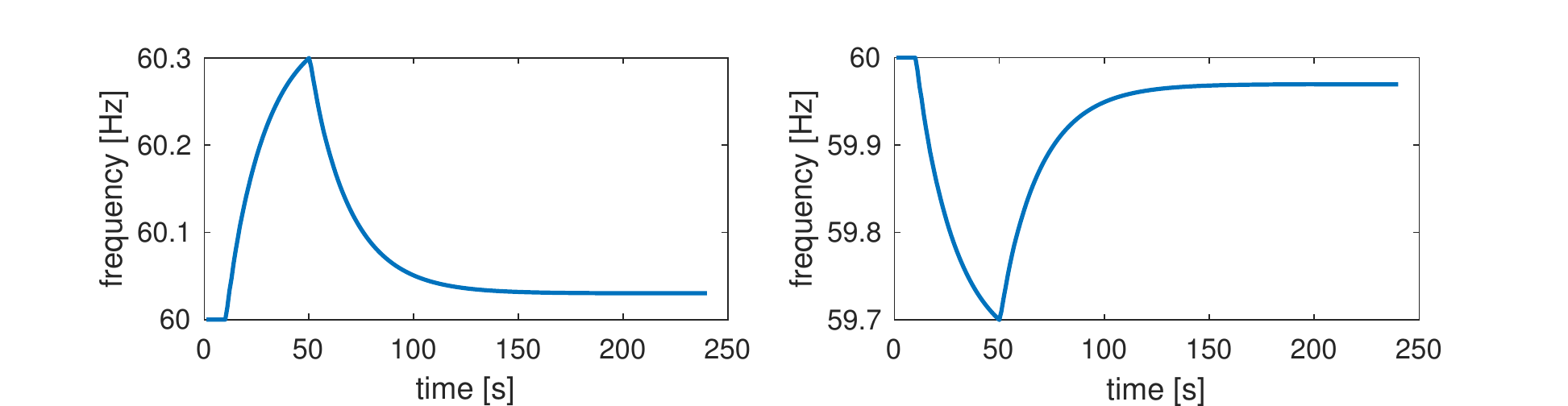}
\caption{Sample under- and over-frequency events generated in IEEE 39-bus.}
\label{F:freq}
\end{figure}

In order to test the prioritized decentralized frequency response algorithm, a random collection of 1000 ACs and 1000 EWHs was generated\footnote{Parameters for the EWHs were taken from \cite{Kundu:2017Assessment}, while the parameter values for the ACs are drawn randomly from the following range of values: $P\in[5.5,6.5]$\,kW,\,$R\in[2,2.4]$\,$^o$F/kW,\,$C\in[3.24,3.96]$\,kW-hr/$^o$F,\,$T_{set}\in[70,74]$$^o$F,\,$T_a\in[80,95]$$^o$F and $\eta=2.5$,\,}. The limits for the frequency response were chosen as:
\begin{itemize}
\item 59.7\,Hz and 59.995\,Hz for under-frequency response, and
\item 60.005\,Hz and 60.3\,Hz for over-frequency response.
\end{itemize} 
Several frequency events were created in an IEEE 39-bus network, by injecting disruptions via changing the loads. Fig.\,\ref{F:freq} illustrates two such events. Performance of the frequency response control algorithm was tested against these various frequency events (shifted in time as appropriate).

\begin{figure}[thpb]
\centering
\hspace{-0.15in}\subfigure[Response of 1000 ACs]{
\includegraphics[scale=0.25]{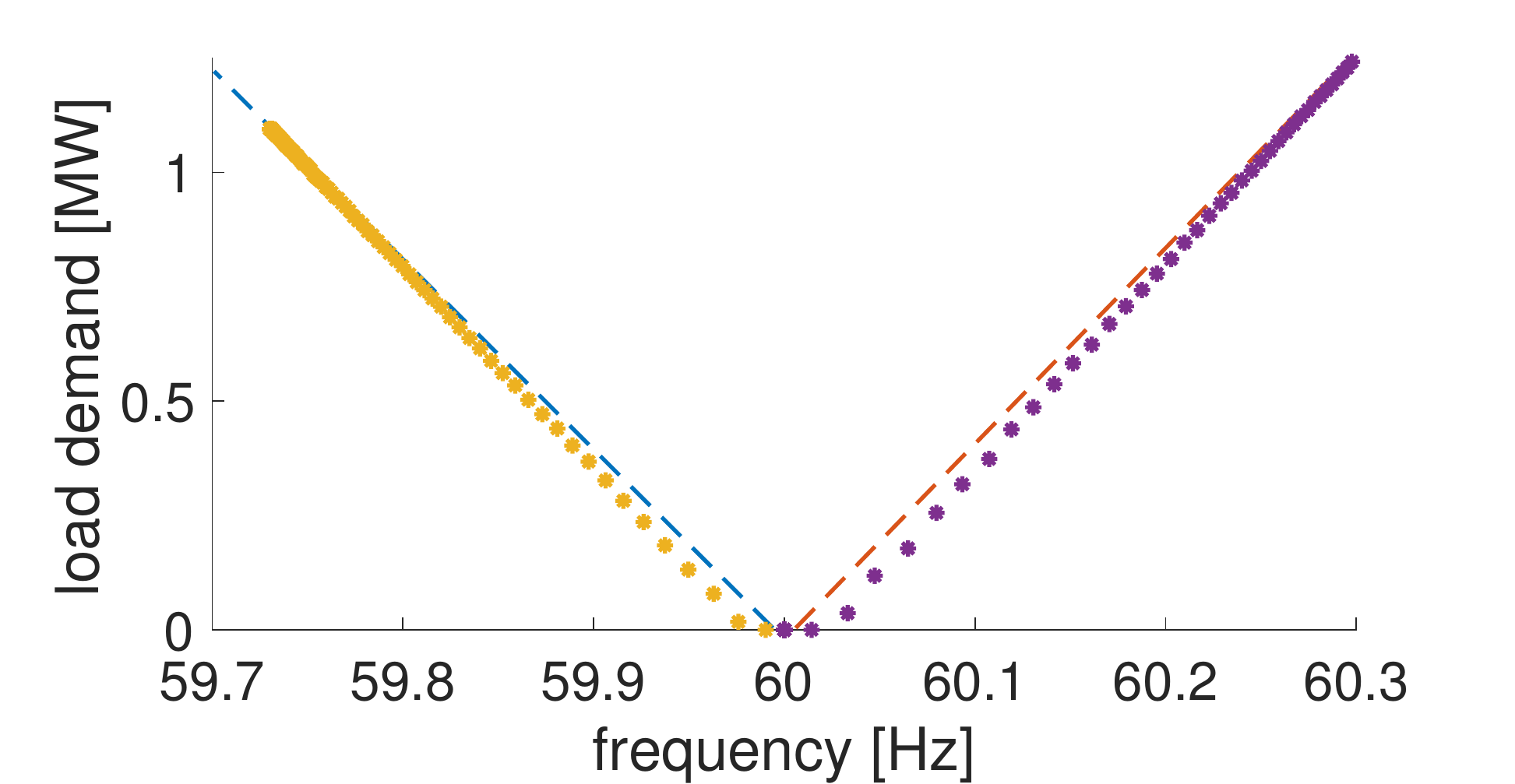}\label{F:respEWH}
}\hspace{-0.3in}
\subfigure[Response of 1000 EWHs]{
\includegraphics[scale=0.25]{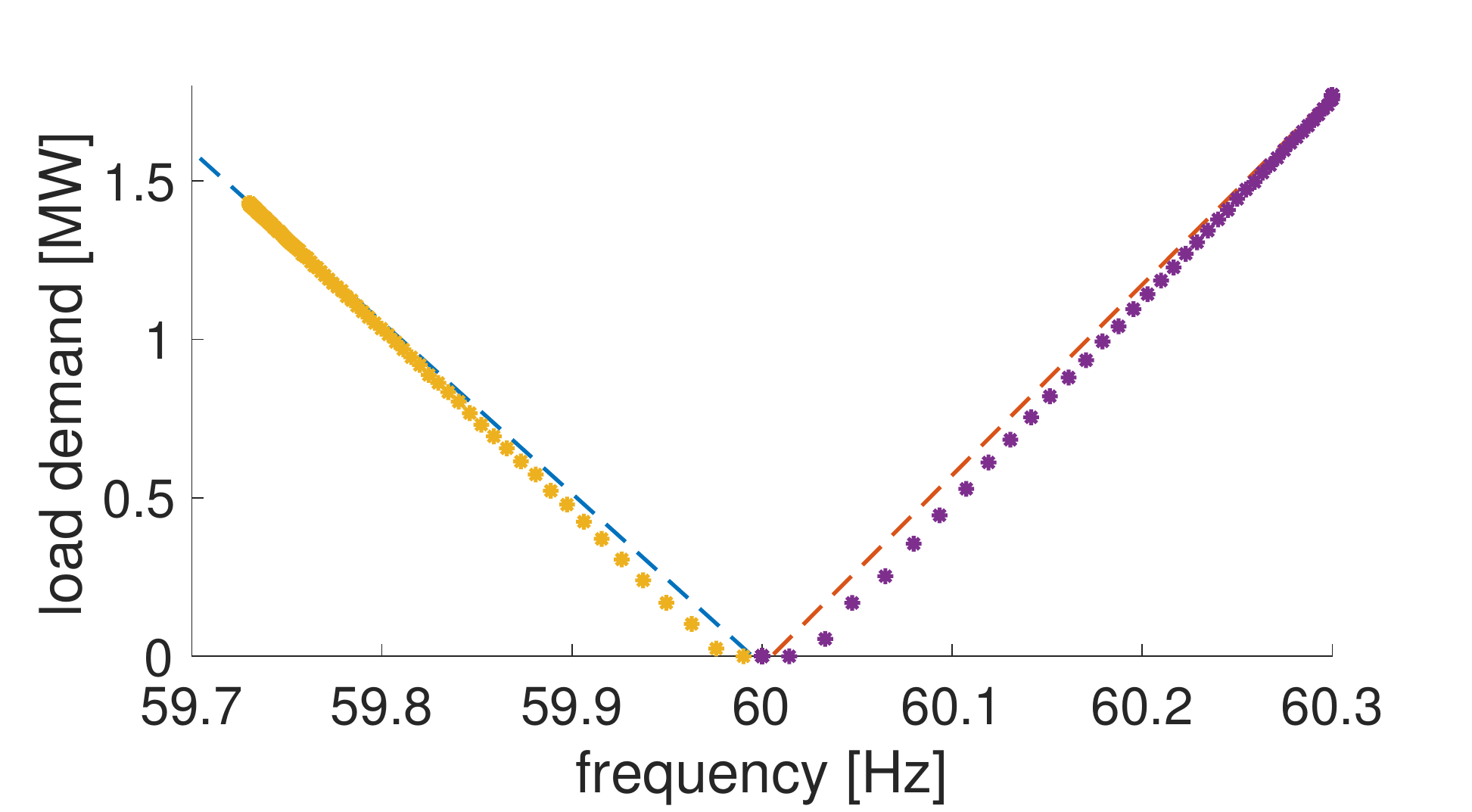}\label{F:respAC}
}\caption[Optional caption for list of figures]{Target (dashed) and achieved (dotted) frequency response curves}
\label{F:resp}
\end{figure}

Fig.\,\ref{F:resp} illustrates how the target and achieved frequency response curves typically look like, simulated using 1000 EWHs (Fig.\,\ref{F:respEWH}) and 1000 ACs (Fig.\,\ref{F:respAC}). The results are obtained for a 5\,min control window, with a sampling time of $\Delta t = 1$\,s, while the response capacity is set at 60\,\% of the maximal capacity as defined in \eqref{E:maxcap}. We note that, with either type of the devices, the achieved response is fairly close to the target response, with some errors being observed for lower frequency deviations. This key observation can be explained by looking at the frequency events in Fig.\,\ref{F:freq}. The rate at which the frequency changes is high when the frequency deviation is low (i.e. closer to the nominal value of 60\,Hz) and gradually decreases (due to damping effect by the generators) as the frequency deviates further. This suggests that the sampling time should be chosen sufficiently small for faster frequency events.

\begin{figure}[thpb]
\centering
\hspace{-0.15in}
\subfigure[Cascading event (left) and response (right).]{
\includegraphics[scale=0.48]{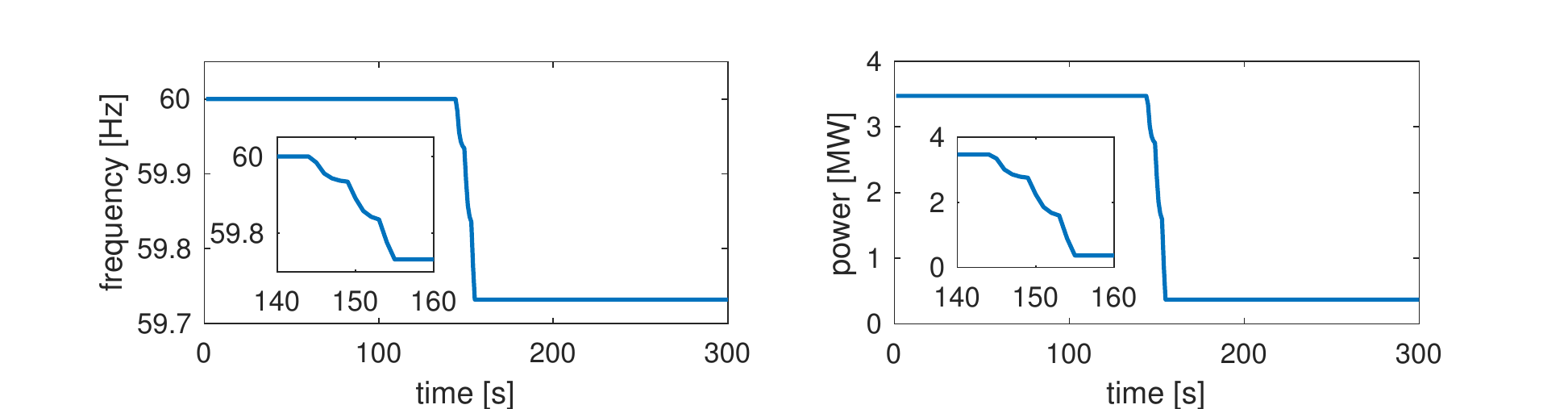}\label{F:resp_profile}
}
\hspace{-0.3in}
\subfigure[Performance metric evaluation.]{
\includegraphics[scale=0.48]{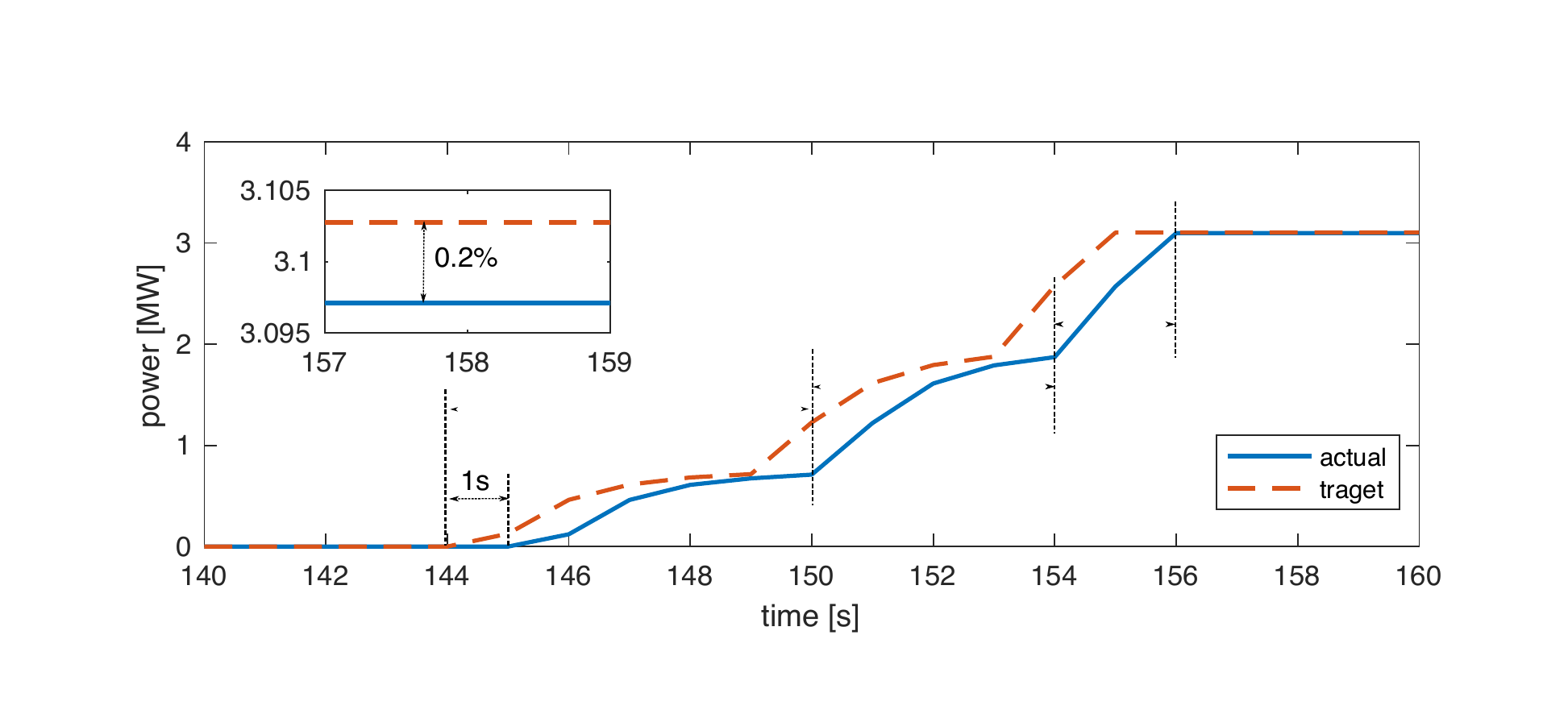}\label{F:resp_perf}
}\caption[Optional caption for list of figures]{Performance under a cascading contingency.}
\label{F:contingency}
\end{figure}

Fig.\,\ref{F:contingency} shows the response of the group of 1000 ACs and 1000 EWHs in response to a cascading contingency where an initial contingency (created by load drop) leads to two subsequent under-frequency events. Control was chosen as 5\,min. The net response magnitude achieved is 3.097\,MW, which is less than 0.2\,\% of the target 3.103\,MW, i.e. RMVT$\leq0.2\%$. Note that, the achieved response is able to track the target response curve always within 1\,s delay. 

\begin{figure}[thpb]
\centering
\includegraphics[scale=0.48]{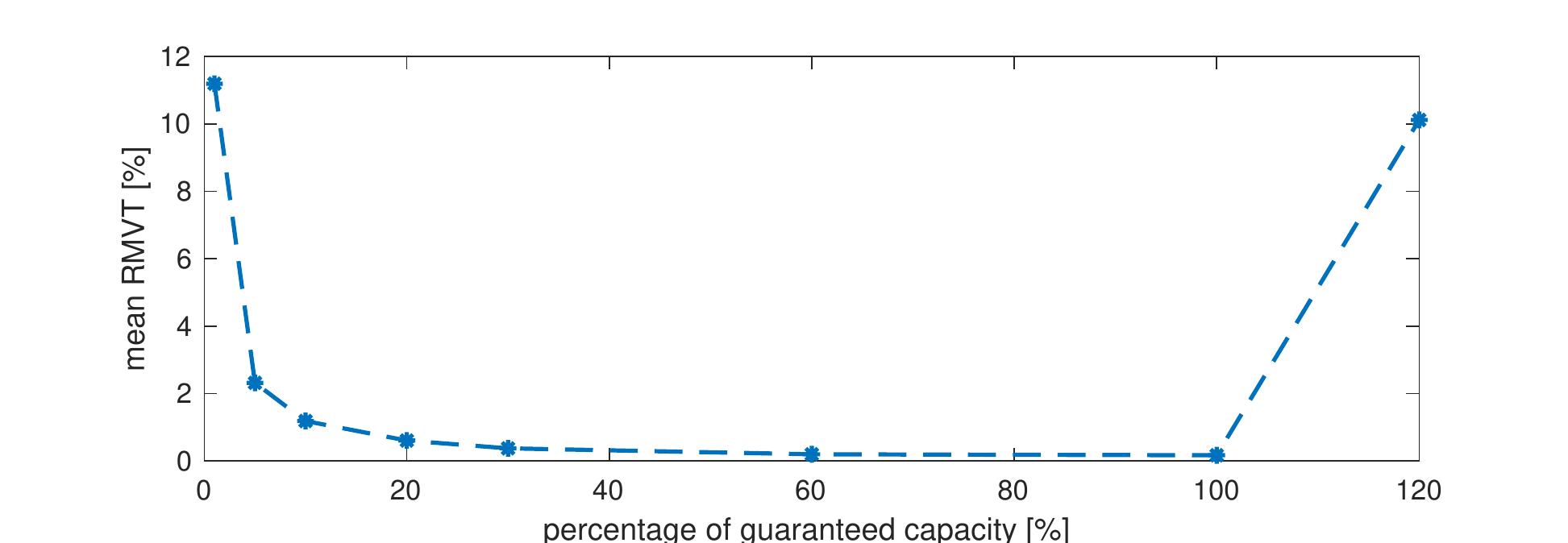}
\caption{Error statistics at varied commitment level as percentage of the guaranteed capacity (computed in \eqref{E:maxcap}).}
\label{F:discrete_err}
\end{figure}

Discrete allocation of frequency thresholds introduce certain error in the response curve. 
Fig.\,\ref{F:discrete_err} shows the error statistics (in the form of mean RMVT) as the committed capacity is varied between 1\,\% to 120\,\% of the maximal guaranteed capacity (given in \eqref{E:maxcap}), i.e. the ratio of $\Delta p^{max}/\left[\Delta p\right]_k^{cap}$ is plotted on the x-axis. Clearly below 20\,\% commitment ($\sim$1\,MW), the error is high due to discrete allocation of resources, which steadies to $<$0.3\,\% (coming from sampling time) from over 20\,\% up to 100\,\%. Beyond 100\,\% commitment, the error increases again due to unavailability of sufficient number of `fit' devices.

Finally Monte-Carlo simulations were run to test the algorithm under various different scenarios. With a fixed (but randomly generated) population of 1000 EWHs and 1000 ACs, different scenarios were created by changing the initial operating condition, the length of the control window as well as the time when the frequency event happens. The RMT level was chosen to be 60\,\% of the estimated maximal guaranteed capacity (computed from \eqref{E:maxcap}). In Table\,\ref{Tab:priority} the mean RMVT values for 6 different scenarios are presented. Across all these scenarios, the priority-based algorithm performs uniformly well with mean RMVT $<0.3\,\%$. However, when not using the priority-based algorithm, the RMVT increases and is sensitive to the length of the control window and when the time of occurrence of the event.

\begin{table}[thpb]
\caption{Mean RMVT [\%] with priority-based allocation}
\label{Tab:priority}
\begin{center}
\begin{tabular}{|l*{4}{|c}}\hline
\backslashbox{Control Window}{Event Time} 
& `start' & `middle' & `end' \\\hline\hline
5\,min & 0.2078 & 0.2020 &  0.2021 \\\hline
15\,min & 0.2437 & 0.2602 & 0.2637\\\hline                             
\end{tabular}
\end{center}
\end{table}

\begin{table}[thpb]
\caption{Mean RMVT [\%] without priority-based allocation}
\label{Tab:nopriority}
\begin{center}
\begin{tabular}{|l*{4}{|c}}\hline
\backslashbox{Control Window}{Event Time} 
& `start' & `middle' & `end' \\\hline\hline
5\,min & 1.5747 & 1.4427  &  1.3613 \\\hline
15\,min & 1.0400 & 1.3960 & 5.7700\\\hline                                                           
\end{tabular}
\end{center}
\end{table}

\section{Conclusion}\label{S:concl}
In this paper we discuss a hierarchical control framework where load aggregators coordinate ensembles of controllable switching loads to commit certain frequency responsive reserve to the grid operator. At the start of a control window, the aggregator assigns a frequency threshold to each switching device which it uses to respond autonomously to frequency events. We present a metric to evaluate the fitness of each device in providing autonomous frequency response, leveraging local device-level measurements and load dynamics. The fitness values were used to prioritize loads for response service, with performance improvement over the non-prioritized method of threshold allocation. Simulation results are provided to validate the control performance under various operational conditions. Future work will extend the deterministic formulation to incorporate stochastic uncertainties in the loads behavior.


\section*{Acknowledgment}

This work was carried out (contract DE-AC02-76RL01830) under the support from the U.S. Department of Energy as part of their ARPA-E/NODES program.

\IEEEtriggeratref{6}


\bibliographystyle{IEEEtran}
\bibliography{IEEEabrv,RefList,RefKundu}
%
%

\end{document}